\documentclass[10pt]{article}
\usepackage{times}

\usepackage{latexsym}
\usepackage{amsmath}
\usepackage{amssymb}

\newtheorem{Theor}{Theorem}[section]
\newtheorem{Prop}[Theor]{Proposition}

\newtheorem{Lemma}[Theor]{Lemma}

\newtheorem{Def}[Theor]{Definition}

\def\Dem{\hbox{\textsc {Proof.}\,\,}}
\def\Aut{\hbox{\sffamily {Aut}\,}} 
\def\Id{\hbox{\sffamily {Id}\,}}
\def\Tri{\hbox{\sffamily {Tri}\,}}

\def\Min{\hbox{\sffamily {Min}\,}} 
\def\rank{\hbox{\rm rank \,}}

\def\Der{\hbox{\sffamily {Der}\,}}

\def\supp{\,\hbox{\rm supp\,}}

\def\ba{\hbox{{\bf a}}}

\def\qed{\hfill $\Box $}

\usepackage{version}

\newenvironment{mm}
  {\marginpar{$\Downarrow$}
}
  {\marginpar{$\Uparrow$}}

\def\nomarks{\def\mm{}\def\endmm{}}

\nomarks

\begin{document}

\title{Manifolds of algebraic elements in the algebra
  $\mathcal{L}(H)$ \\ of bounded linear operators.}

\author{Jos\'e M. Isidro \thanks{Supported by Ministerio de 
Educaci\'on y Cultura of Spain, Research Project PB 98-1371.}\\
Facultad de Matem\'aticas, \\
Universidad de Santiago,\\
Santiago de Compostela, Spain.\\
{\normalsize\tt jmisidro@zmat.usc.es}
}

\date{October 10,  2001}

\maketitle

\begin{abstract}
  Given a complex Hilbert space $H$, we study the differential
  geometry of the manifold $\mathcal{A}$ of 
		normal algebraic elements 
  in $Z=\mathcal{L}(H)$. We represent $\mathcal{A}$ as a disjoint
	 union of connected subsets $\mathcal{M}$ of $Z$. Using the 
		algebraic structure of $Z$, a 
  torsionfree affine connection $\nabla$ (that is invariant under the
  group $\Aut (Z)$ of automorphisms of $Z$) is defined on each of these
		connected components and the geodesics are computed. 
  In case $\mathcal{M}$ consists of elements that have a
		fixed finite rank $r$, $(0<r<\infty)$, $\Aut (Z)$-invariant Riemann
		and K\"ahler structures are defined on $\mathcal{M}$ which in this way	becomes 
  a totally geodesic symmetric holomorphic manifold. 
	 \bigskip

\noindent {\bf Keywords.} 
JB$^*$-triples, Grassmann manifolds, Riemann manifolds.\\

\noindent {\bf AMS 2000 Subject Classification.} 48G20,  72H51. 
\end{abstract}


\section{Introduction}

In this paper we are concerned with the differential geometry of some  
infinite-dimensional Grassmann manifolds  
in $Z\colon = 
\mathcal{L}(H)$, the space of 
bounded linear operators $z\colon H\to H$ in a complex 
Hilbert space $H$. Grassmann manifolds are a classical object in 
Differential Geometry and in recent years several authors have considered 
them in the Banach space setting. Besides the Grassmann structure, a 
Riemann and a K\"ahler structure has sometimes been defined even in the 
infinite-dimensional setting. Let us recall some aspects of the 
topic that are relevant for our purpose. 

The study of the manifold of minimal projections in a finite-dimensional simple 
formally real Jordan algebra was made by U. Hirzebruch in \cite{HIR}, who 
proved that such a manifold is a compact symmetric Riemann space of rank 1, and 
that every such a space arises in this way. Later on, Nomura in \cite{NOM1, NOM2} 
established similar results for the manifold of fixed finite rank projections in a 
topologically simple real Jordan-Hilbert algebra. In \cite{ISMA}, the authors 
studied the Riemann and K\"ahler structure of the manifold of finite rank projections in
$Z$ without the use of any global scalar product. As pointed out there, the Jordan-Banach 
structure of $Z$ encodes information about the differential geometry of some 
manifolds naturally associated to it, one of which is the manifold of algebraic elements 
in $Z$. On the other hand, the 
Grassmann manifold of all projections in $Z$  
has been discussed by Kaup in \cite{KAUG} and \cite{KAUJ}. 
See also \cite{CHIS,ISI} for related results.

It is therefore reasonable to ask whether a Riemann structure  can  
be defined in the set of algebraic elements in $Z$, and how
does it behave when it exists. We restrict our considerations to the set $\mathcal{A}$ 
of all normal algebraic elements in $Z$ that have finite rank. Remark that the assumption 
concerning the finiteness of the rank can not be dropped, as proved in \cite{ISMA}. 
Normality allows us to use spectral 
theory which is an essential tool. In the case $H=\mathbb{C}^n$, all elements in 
$Z$ are algebraic (as any square matrix is a root of its characteristic polynomial) 
and have finite rank. 
Under the above restrictions $\mathcal{A}$ is represented as a 
disjoint union of connected subsets $M$ of $Z$, each of which is invariant
under $\Aut (Z)$ (the group of all  C$^*$-automorphisms of $Z$). Using algebraic
tools, a holomorphic manifold structure and an $\Aut (Z)$-invariant affine connection $\nabla$
are introduced on $M$ and its geodesics are calculated. One of the 
novelties is that we take JB$^*$-triple system approach instead of the Jordan-algebra approach
of 
\cite{NOM1, NOM2}. As noted in  \cite{CHIS} and \cite{ISI}, within this 
context the algebraic structure of JB$^*$-triple acts as a substitute for the 
Jordan algebra structure. In case $M$ consists of elements that 
have a fixed finite rank $r$, $(0<r<\infty)$, the JB$^*$-triple structure provides a {\sl local
scalar product} known as the {\sl algebraic metric} of Harris (\cite{DIN}, prop. 9.12).
Although $Z$ is not a Hilbert space, the use of the algebraic scalar product allows us to
define an $\Aut (Z)$-invariant Riemann and a K\"ahler structure on $M$. We
prove that
$\nabla$ is the Levi-Civita and the K\"ahler connection of $M$, and that $M$ is a symmetric
holomorphic manifold on which $\Aut ^{\circ}(Z)$ acts transitively as a group 
of isometries. 

The role that projections play in the study of the 
algebra $Z=\mathcal{L}(H)$ is taken by tripotents in the study
of a JB$^*$-triple system. A spectral calculus and a
notion of algebraic element is available in
the stetting of JB$^*$-triples, and the manifold of all 
finite rank algebraic elements in a JB$^*$-triple $Z$ is studied in 
the final section.


\section{Algebraic preliminaries.}

For a complex Banach space $X$ denote by $X_{\mathbb{R}}$ the underlying
real Banach  space, and let $\mathcal{L}(X)$ and
$\mathcal{L}_{\mathbb{R}}(X)$ respectively be the Banach algebra of
all bounded complex-linear operators on
$X$ and the Banach algebra of all bounded real-linear
operators on $X_{\mathbb{R}}$. A complex Banach space $Z$ with a 
continuous mapping $(a, b, c)  
\mapsto \{a b c\}$ from $Z\times Z\times Z$ to $Z$ is called a {\it JB*-triple} 
if the following conditions are satisfied for all $a, b, c, d \in Z$, where 
the operator $a\square b\in \mathcal{L}(Z)$ is defined by $z\mapsto \{abz\}$ and 
$\lbrack\, , \, \rbrack$ is the commutator product:
\begin{enumerate}
\item $\{abc\}$ is symmetric complex linear in $a, c$ and conjugate linear 
in $b$.
\item $\lbrack a\square b , \, c\square d \rbrack = \{abc\}\square d -
c\square \{dab\}.$
\item $a\square a$ is hermitian and has spectrum $\geq 0.$
\item $\Vert \{aaa\}\Vert = \Vert a\Vert ^3$. 
\end{enumerate}
If a complex vector space $Z$ admits a JB*-triple structure, then the norm and 
the triple product determine each other. For $,x,y, z\in Z$ we write $L(x,y)(z) = 
(x\Box y)(z)$ and $Q(x,y)(z)\colon =\{xzy\}$. Note that $L(x,y)\in \mathcal{L}(Z)$
whereas $Q(x,y)\in \mathcal{L}_{\mathbb{R}}(Z)$, and that the operators 
$L_a= L(a,a)$ and $Q_a=Q(a,a)$ commute. 
A {\it derivation} of a JB*-triple $Z$ is an element $\delta \in \mathcal{L}(Z)$ such
that 
$
\delta \{z z z \}= \{(\delta z)  z  z\} + \{z  (\delta z) z\} +
\{z z (\delta z)\} \; 
$ 
and an {\it automorphism} is a bijection $\phi \in \mathcal{L}(Z)$ such 
that $\phi \{zz z\}= \{(\phi z) ( \phi z)(\phi z)\}$ for $z\in Z$. The latter
occurs if and only if $\phi$ is a surjective linear isometry of $Z$. 
The group $\Aut (Z)$ of automorphisms of
$Z$ is a real Banach-Lie 
group whose Banach-Lie algebra is the set $\Der(Z)$ of all derivations of $Z$. 
The connected component of the identity in $\Aut (Z)$ is denoted by 
$\Aut^{\circ}(Z)$. 
Two elements $x, y \in Z$ are {\sl orthogonal} if $x \square y = 0$ 
and $e\in Z$ is called a {\it tripotent} if $\{e e e \}=e$, the set of which is
denoted by $\Tri(Z)$. 
For $e\in \Tri(Z)$, the set of eigenvalues of $e\Box e\in
 \mathcal{L}(Z)$ is contained in 
$\{0,\, {1\over 2}, \, 1\}$ and the topological direct sum 
decomposition, called the {\it Peirce decomposition} of $Z$,  
\begin{equation}\label{pd}
Z=Z_1(e)\oplus Z_{1/2}(e)\oplus Z_0(e).
\end{equation}
holds. Here $Z_k(e)$ is the $k$-
eigenspace of $e\Box e$ and the {\sl Peirce projections} are 
$$
P_1(e) = Q^2(e), \qquad P_{1/2}(e)= 2(e\Box e- Q^2(e)), \qquad 
P_0(e)= \hbox{\rm Id} -2e\Box e+Q^2(e).
$$ 
We will use the {\sl Peirce rules}
$ \{Z_i (e)\, Z_j (e)\, Z_k (e) \} \subset Z_{i-j+k} (e) $
where $Z_l (e) = \{0\}$ for $l \neq 0,1/2,1$. In particular, 
every Peirce space is a JB$^*$-subtriple of $Z$ and $Z_1(e)\Box Z_0(e)=\{0\}$. 
We note that $Z_1(e)$ is a complex unital JB*-algebra in 
the product $a\circ b\colon =
\{ a eb\}$ and involution $a^{\#} \colon =\{eae\}$. 
 Let
$$
A(e)\colon = \{ z\in Z_1(e)\;\colon \; z^{\#}=z\}.
$$
Then we have $Z_1(e)=A(e)\oplus iA(e)$. 
A tripotent $e$ in a JB$^*$-triple $Z$ is said to be {\sl minimal} if
$e\neq 0$ and 
$ P_1(e)Z=\mathbb{C}\, e$, and we let $\Min (Z)$ be the set of them. 
If $e\in \Min(Z)$ then $\Vert e\Vert=1$. 
A JB$^*$-triple $Z$ may have no non-zero tripotents.

Let ${\bf e}=(e_1,\cdots ,e_n)$ be a finite sequence of non-zero mutually 
orthogonal tripotents $e_j\in Z$, and define for all integers $0\leq
j,k\leq n$ the linear subspaces
\begin{equation}\label{jpd}
\begin{split} 
Z_{j,j}({\bf e})&= Z_1(e_j)\qquad \qquad \qquad \quad \;\;1\leq j\leq n,\\
Z_{j,k}({\bf e})= Z_{k,j}({\bf e})&=Z_{1/2}(e_j)\cap Z_{1/2}(e_k) 
\qquad 1\leq j, k \leq n, \;\; j\neq k, \\
Z_{0,j}({\bf e})=Z_{j,0}({\bf e})&= Z_1(e_j)\cap \bigcap_{k\neq j}
Z_0(e_k)\qquad 1\leq j\leq n,\\
Z_{0,0}({\bf e})&= \bigcap_jZ_0(e_j).
\end{split}
\end{equation}
Then the following topologically direct sum decomposition, called the 
Peirce decomposition relative to ${\bf e}$, holds
\begin{equation}\label{jds}
 Z= \bigoplus_{0\leq j \leq k\leq n}Z_{j,k}({\bf e}).
\end{equation}
The Peirce spaces multiply according to the rules
$\{Z_{j,m}Z_{m,n}Z_{n,k}\}\subset Z_{j,k}$, and all products that cannot
be brought to this form (after reflecting pairs of indices if necessary) 
vanish. 
In terms of this decomposition, the Peirce spaces of the tripotent
$e \colon = e_1+\cdots +e_n$ are 
\begin{equation}
\begin{split}
Z_1(e )= \bigoplus_{j,k}Z_{j,k}({\bf e})&= \Big(\bigoplus_{1\leq
j\leq n}Z_{j,j}({\bf e})\Big) \; \oplus \;
\Big(\bigoplus_{\substack{1\leq j, k\leq n\\ j\neq k}}Z_{j,k}({\bf e}) 
\Big), \\
Z_{1/2}(e )&=\bigoplus_{1\leq j\leq n}Z_{0,j}({\bf e}), \qquad \quad 
Z_0({\bf e})=Z_{0,0}({\bf e}).
\end{split}
\end{equation}

Recall that every C*-algebra $Z$ is a JB*-triple with respect
to the triple product $2\{abc\} \colon =(ab^*c+cb^*a)$. In that case, 
every  projection in $Z$ is a tripotent and more generally the tripotents
are precisely the partial isometries in $Z$. C$^*$-algebra derivations
and C$^*$-automorphisms are derivations and automorphisms of $Z$ as a 
JB$^*$-triple though the converse is not true. 

We refer to \cite{KAUR}, \cite{KAUJ}, \cite{LOOS}, 
\cite{UPM} and the references therein for
the background of JB$^*$-triples theory.


\section{Manifolds of algebraic elements in $\mathcal{L}(H).$} 

From now on, $Z$ will denote the C$^*$-algebra $\mathcal{L}(H)$. An element $a\in Z$ 
is said to be {\sl algebraic} if it satisfies  the equation $p(a)=0$ for
some non identically null polynomial 
$p\in \mathbb{C}[X]$. By elementary spectral theory $\sigma (a)$, the spectrum of $a$ in
$Z$, is a finite set whose elements are roots of the algebraic 
equation 
$p(\lambda )=0$. In case $a$ is {\sl normal} we have 
\begin{equation}\label{ae}
a= \sum_{\lambda \in \sigma (a)}\lambda \,e_\lambda
\end{equation}
where $\lambda$ and $e_\lambda$ are, respectively, the spectral values
and the corresponding spectral projections of $a$. If $0\in \sigma (a)$ then $e_0$,
the projection onto $\ker (a)$, satisfies $e_0\neq 0$ but in (\ref{ae}) 
the summand $0\,e_0$ is null and will be omitted. In
particular, in (\ref{ae}) the numbers $\lambda$ are non-zero pairwise
distinct complex numbers and the 
$e_\lambda$ are pairwise orthogonal non-zero projections. 
We say that
$a$ has {\sl finite rank} if $\dim a(H)<\infty$, which always occurs if $\dim (H)<\infty$.
Set $r_\lambda \colon = \rank (e_\lambda )$. Then $a$ has finite rank if and only if 
$r_\lambda <\infty$ for all $\lambda\in\sigma (a)\backslash \{0\}$ (the case
$0\in \sigma (a)$ and $\dim \ker a=\infty$  may occur).

Thus, every finite rank normal algebraic element $a\in Z$ gives rise to:  
(i) a positive integer $n$ which is the cardinal of $\sigma (a)\backslash \{0\}$, (ii) an
ordered  n-uple $(\lambda_1,\cdots ,\lambda_n)$ of numbers in $\mathbb{C}\backslash \{0\}$
which is the set of the pairwise distinct non-zero spectral values of $a$,
(iii) an ordered n-uple
$(e_1,\cdots , e_n)$ of non-zero pairwise orthogonal projections, and (iii) an ordered 
n-uple $(r_1,\cdots , r_n)$ where $r_k\in \mathbb{N}\backslash \{0\}$.

The spectral resolution of $a$ is
unique except for the order of the summands in (\ref{ae}), therefore these
three n-uples are uniquely determined up to a permutation of the indices $(1,\cdots ,n)$.
The operator $a$ can be
recovered from the set of the first two ordered n-uples,
$a$ being given by (\ref{ae}). 

Given the n-uples $\Lambda\colon =(\lambda_1,\cdots , \lambda_n)$ and $R\colon = 
(r_1,\cdots ,r_n)$ in the above conditions, we let
\begin{equation}\label{nn}
M(n,\,\Lambda ,\, R)\colon = 
\{\,\sum_k\lambda_ke_k\;\colon \;
\; e_je_k=0 \;\; \hbox{\rm for}\; j\neq k,\;  
\rank (e_k)=r_k, \; 1\leq j,\, k\leq n \,\} 
\end{equation}
be the set of the elements (\ref{ae}) where the coefficients $\lambda_k$ and ranks $r_k$
are given and the $e_k$ range over non-zero, pairwise orthogonal projections of rank 
$r_k$. For instance, for $n=1$, $\Lambda =\{1\}$ and
$R=\{r\}$ we obtain the manifold of projections with a given finite rank $r$, that
was studied  in \cite{ISMA}. For the n-uple $\Lambda =(\lambda_1,\cdots , \lambda_n)$ 
we set $\Lambda^*\colon =(\bar \lambda_1, \cdots , \bar \lambda_n)$. The involution 
$z\mapsto z^*$ on $Z$ induces a map $M(n, \Lambda , R)\to M(n, \Lambda^* , R)$ where 
$M(n, \Lambda , R)^* = \{z^*\colon z\in M\}=M(n, \Lambda^* , R)$, and $\Lambda\subset 
\mathbb{R}$ if and only if $M(n, \Lambda , R)$ consists of hermitian elements.  

For a normal algebraic element $a=\sum_{\lambda\in \sigma (a)\backslash \{0\}}\lambda
e_\lambda$ we define its {\sl support} to be  the projection  
$$ 
{\bf a} =\supp a \colon =\sum_{\lambda\in \sigma(a)\backslash \{0\}}e_\lambda = e_1+\cdots
+e_n.
$$
It is clear that $h (\supp (a)) = \supp h(a)$ holds for all $h\in \Aut^{\circ}(Z)$,
which  combined with the $\Aut^{\circ}(Z)$-invariance of Peirce projectors $P_k$ gives the
following useful formula
\begin{equation}\label{spi}
P_k\big(\supp h(a)\big)=P_k\big(h \supp (a)\big) =h\,P_k\big(\supp (a)\big)\,h^{-1}, \qquad
(k= 1, 1/2, 0).
\end{equation}
\begin{Prop}\label{al}
Let $\mathcal{A}$ be the set of all normal algebraic elements of finite rank in $Z$, 
and let $M(n,\Lambda , R)$ be defined as in (\ref{nn}). Then 
\begin{equation}\label{eq} 
\mathcal{A}=\bigcup_{n,\,\Lambda ,\,R}M(n,\,\Lambda ,\, R)
\end{equation}
is a disjoint union of $\Aut ^{\circ}(Z)$-invariant connected subset of $Z$  
on which  
the group $\Aut ^{\circ}(Z)$ acts transitively. 
\end{Prop}
\Dem 
We have seen before that $\mathcal{A}\subset\bigcup_{n,\,\Lambda ,\,R}M(n,\,
\Lambda ,\, R)$. Conversely, let $a$ belong to some $M(n,\,\Lambda ,\, R)$ hence 
we have $a= \sum_k\lambda_ke_k$ for some orthogonal projections $e_k$. Then $\Id =
(e_1+\cdots +e_n) +f$ where $f$ is the projection onto $\ker (a)$ in case $0\in
\sigma (a)$ and $f=0$ otherwise. The above properties of the 
$e_k, f$ yield easily $ap(a)=0$ or $p(a)=0$ according to the cases, where $p\in
\mathbb{R}[X]$ is the polynomial $p(z)= (z-\lambda_1). \cdots .(z-\lambda_n)$. Hence $a\in 
\mathcal{A}$. Clearly (\ref{eq}) is union of disjoint subsets.

Fix one of the sets $M\colon = M(n,\,\Lambda ,\, R)$ and take any pair $a, b\in M$. 
Then 
$$ a= \lambda_1p_1+\cdots +\lambda_np_n, \qquad b= \lambda_1q_1+\cdots +\lambda_nq_n.$$
In case $0\in \sigma (a)$, set $p_0\colon =\Id -\sum_kp_k$ and $q_0\colon =\Id
-\sum_kq_k$. Since
$\rank p_k=\rank q_k$, the projections $p_k$ and $q_k$ are unitarily equivalent and so 
are $p_0$ and $q_0$. Let us choose orthonormal basis
$\mathcal{B}_k^p$ and $\mathcal{B}_k^q$ in the ranges $p_k(H)$ and $q_k(H)$ for $k=
0,\,1,\cdots ,n$. Then 
$\bigcup_k \mathcal{B}_k^p$ and $\bigcup_k \mathcal{B}_k^q$ are two orthonormal basis in
$H$. The unitary operator $U\in Z$ that exchanges these basis satisfies $Ua=b$. In
particular 
$M$ is the orbit of any of its points under the action of the unitary group of $H$. 
Since this group is connected and its action on $Z$ is continuous, $M$ is connected. 
\qed

Let $a\in Z$ be a normal algebraic element with finite rank and  
${\bf a}=\supp (a)$ its support. 
In the Peirce decomposition 
\begin{equation*}
Z= Z_1({\bf a})\oplus Z_{1/2}({\bf a})\oplus Z_0({\bf a})
\end{equation*}
every Peirce space $Z_k({\bf a})_s$ is invariant under the natural involution 
$^*$ of $Z$, and we let $Z_k({\bf a})_s$ denote its selfadjoint part,  
$(k=1,1/2,0)$.  
In what follows, the map $Z\times Z\to Z$ given by $(x,y)\mapsto 
g(\ba ,x)y$, and the partial maps obtained by fixing 
one of the variables, will play an important role. 
For every fixed value $x\in Z_{1/2}({\bf a})$, we get an operator 
$g(\ba ,x)(\cdot )$ which is an inner JB$^*$-triple derivation of 
$Z$, hence we have an operator-valued continuous real-linear map
$Z_{1/2}({\bf a}) \to \Der (Z)$. Moreover $g(\ba ,x)(\cdot )$ is 
a C$^*$-algebra derivation if and only if $x\in Z_{1/2}({\bf a})_s$ 
(see \ref{der}). 
For $y=a$ fixed, we get the map 
$x\mapsto g(\ba ,x)a$ for which we introduce the notation
$$
\Phi_a (x)\colon = g(\ba ,x)a= \{\ba \,x\,a\}-\{x\,\ba \,a\}= \big( Q(\ba ,a)
-L(\ba ,a)\big) x, \qquad x\in Z.
$$
First we discuss $Z_{1/2}({\bf a})$.

\begin{Prop}\label{sp}
Let $a\in Z$ be a normal algebraic element of finite rank, and let 
${\bf a}= e_1+\cdots +e_n$ be its support. Then 
$Z_{1/2}({\bf a})$ consists of the  operators 
\begin{equation}\label{1/2}
u=  \sum_k u_k, \quad  u_k\in Z_{1/2}(e_k), 
\quad  e_k\,u_j=u_j\,e_k=0, \quad j\neq k, \quad (1\leq j,\,k\leq n).
\end{equation}
If $u\in Z_{1/2}({\bf a})_s$, then we have the additional condition $u_k\in 
Z_{1/2}(e_k)_s$. 
\end{Prop} 
\Dem Let $u\in Z$ be selfadjoint. The relation $u\in Z_{1/2}({\bf a})$ is 
equivalent to $u=2\{{\bf a} {\bf a} u \}$ which now reads 
$$
u=2\{{\bf a} {\bf a}u \}={\bf a}{\bf a}^*u+u{\bf a}^*{\bf a}= 
\sum_k(e_ku+ue_k)= \sum _ku_k
$$ 
where 
\begin{equation}\label{df}
u_k\colon = e_ku+ue_k \quad \hbox{\rm for} \quad 1\leq k\leq n.
\end{equation} 
Note that $e_j, e_k\in 
Z_1({\bf a})$, hence by the Peirce
multiplication rules
$\{e_jue_k\}\in \{ Z_1({\bf a})Z_{1/2}({\bf a})Z_1({\bf a}) \}= 
\{ 0\}$, that is 
$ e_jue_k+e_kue_j=0 $ for all $1\leq j,\, k \leq n.$  
Multiplying the latter by $e_j$ with $j\neq k$ yields 
$e_jue_k=0$ for $j\neq k,\; (1\leq j,\, k \leq n).$ 
Therefore by (\ref{df}),  
\begin{eqnarray*} 
2\{e_ke_ku_k\}=e_k (e_ku+ue_k)&+(e_ku+ue_k)e_k=\\
(e_ku+ue_k)+2e_kue_k&= 
(e_ku+ue_k)=u_k
\end{eqnarray*}
which shows $u_k\in Z_{1/2}(e_k)$ and clearly $u_k=u_k^*$ for $1\leq k\leq n.$ 
Multiplying in (\ref{df}) by $e_j$ with $j\neq k$ we get $u_ke_j=e_ju_k=0$ 
and in particular $e_j\Box u_k=u_k\Box e_j=0$ for $j\neq k.$  

Conversely, let $u_k$ satisfy the properties in (\ref{1/2}). Then 
$u\colon =\sum_ku_k$ is selfadjoint and 
$e_ku=e_k\, (\sum_j u_j)= e_ku_k$. Similarly$ ue_k= u_ke_k$, hence 
$2\{ {\bf a}{\bf a}u \}={\bf a}{\bf a}^*u+u{\bf a}^*{\bf a}= 
(\sum_je_j)\,u+u\,(\sum_je_j)= \sum_j(e_ju+ue_j)= u,$. 

Using the $^*$-invariance of $Z_{1/2}(\bf a)$ every element in this 
space can be written in the form $u=u_1+iu_2$ with $u_1,u_2\in Z_{1/2}(\bf a)_s$ 
and the result follows easily. 
\qed

The following result should be compared with (\cite{CHIS}, th. 3.1)
\begin{Prop}\label{der}
Let $a\in Z$ be a normal algebraic element of finite
rank and ${\bf a}\colon =\supp(a)$. Then for any  
$u\in Z_{1/2}({\bf a})$, the operator $g({\bf a}, u) \colon = {\bf a}\Box u-
u\Box {\bf a}$ is an inner C$^*$-derivation of $Z$ if and only if 
$u$ is selfadjoint.  
\end{Prop}
\Dem 
Let $a= \sum_k\lambda_k e_k$ and 
${\bf a}= \sum_k e_k$ be the spectral resolution and the support of $a$. 
Suppose $u=u^*$. By 
(\ref{sp}) $u$  has the form $u=\sum u_k$ with $u_k\in Z_{1/2}(e_k)_s$ and 
$ e_k\Box u_j= u_j\Box e_k=0$ for all $j\neq k$. 
Therefore 
\begin{equation}\label{g}
g({\bf a},u)=\sum_k  (e_k\Box u_k-u_k\Box e_k)= \sum_k g(e_k,u_k).
\end{equation}
Here the $e_k$ are projections in $Z$ and $u_k\in Z_{1/2}(e_k)_s$, 
hence by (\cite{CHIS}, th. 3.1) each $g(e_k,u_k)$ is an inner
C$^*$-derivation of $Z$ and so is the sum. Conversely, since $\bf a$ 
is a projection, whenever $g({\bf a},u)$ is a C$^*$-algebra derivation 
we have $u\in Z_{1/2}(\bf a)_s$ again by (\cite{CHIS}, th. 3.1). \qed 
\par 
Now consider the joint Peirce decomposition of $Z$ relative to the 
family $(e_1, \cdots , e_n)$ where $a= \lambda_1e_1+ \cdots +\lambda_ne_n$ 
is the spectral resolution of $a$. Remark that 
$\bigoplus_{1\leq k\leq n}i\,A(e_k)\subset Z_1({\bf a})$ is a direct summand 
of $Z$, hence so is the space 
$$X\colon =  \big(\bigoplus_{1\leq k\leq n}i\,A(e_k)\big) \oplus Z_{1/2}({\bf a}).$$ 
\begin{Prop}\label{lh}
Let $a\in Z$ be a normal algebraic element of finite
rank and ${\bf a}\colon =\supp(a)$. Then $\Phi_a$ is a surjective complex linear 
homeomorphism of $Z_{1/2}({\bf a})$. If $a$ is hermitian then $\Phi_a$ is a  
surjective real linear homeomorphism of $X$ that preserves the subspace 
$\bigoplus_{1\leq k\leq n}i\,A(e_k)$.  
\end{Prop}
\Dem Let $x= iv+u\in X$ where $v\in \bigoplus_{1\leq k\leq n}A(e_k)$ 
and $u\in Z_{1/2}({\bf a})$. The Peirce multiplication rules give for $v= 
\sum_j v_j$ with $v_j\in A(e_j)$ and $u= \sum_k u_k$ according to (\ref{sp}) 
\begin{align}
{}&\{\ba Z_{1/2}(\ba )a\} =\{0\}, \nonumber \\
{}&\{\ba \,iv \,a\}= -i \{ \sum_j e_j \sum_k v_k \sum_l \lambda_l e_l\}= -i \sum_k 
\lambda_k v_k, \nonumber \\
{}&\{ u \,\ba \,a\} = i \{\sum_j u_j \sum_k e_k \sum_l \lambda_l e_l\}= \frac{i}{2} \sum_k 
\lambda_k u_k .\nonumber
\end{align}
Therefore 
\begin{equation}\label{for}
\Phi_a(x) = -2i \sum_k \lambda_kv_k-\frac{1}{2}\sum_k \lambda_ku_k\in 
\big( \bigoplus_{1\leq k\leq n}Z_1(e_k)\big) \oplus  Z_{1/2}(\ba ).  
\end{equation}
It is now clear that $\Phi_a$ preserves $Z_{1/2}({\bf a})$. If $a$ is hermitian then 
$\Lambda \subset\mathbb{R}^n$ and $\Phi_a$ also preserves $\bigoplus_{1\leq k\leq
n}i\,A(e_k)$. Moreover $\Phi_a(x)=0$ with $x\in X$ is
equivalent to 
$\sum \lambda_kv_k=0=\sum \lambda_ku_k$ which is equivalent to 
$v=0=u$ since the coefficients satisfy $\lambda_k\in \sigma (a)
\backslash \{0\}$. We can recover $x$ from 
$\Phi_a(x)$, hence the result follows. \qed 

Recall that a subset $M\subset Z$ is called a {\sl real analytic} (respectively, {\sl
holomorphic}) submanifold  
if to every $a\in M$ there are open subsets $P,Q\subset Z$ and a closed real-linear 
(resp. complex) subspace $X\subset Z$ with $a\in P$ and $\phi (P\cap M)=Q\cap X$
for some bianalytic (resp. biholomorphic) map $\phi \colon P\to Q$. If to every $a\in M$
the linear subspace $X=T_aM$, called the {\sl tangent space} to $M$ at $a$, can
be chosen to be topologically complemented in $Z$ then 
$M$ is called a {\sl direct submanifold} of $Z$. 

Fix one of the sets $M=M(n,\Lambda ,R)$ and a point $a\in M$ with spectral 
resolution $a= \sum_k\lambda_ke_k$. By the orthogonality properties
of the $e_k$, the successive  powers of $a$ have the expression   
$$ a^l= \lambda_1^le_1+\cdots + \lambda_n^le_n, \qquad 1\leq l\leq n,$$
where the determinant $\det (\lambda_k^l)\neq 0$ does not vanish since 
it is a Vandermonde determinant and the $\lambda_k$ are pairwise distinct. 
Thus the $e_k$ are polynomials in $a$ whose   
coefficients are rational functions of the $\lambda_k$.  
Suppose $M$ is a differentiable manifold, and let us obtain its tangent space
$T_aM$. Consider a smooth curve $t\mapsto a(t)$ through $a\in M$, $t\in I$, for   
a neighbourhood $I$ of $0\in \mathbb{R}$ and $a(0)=a$. Each $a(t)$ has a
spectral resolution 
$$
a(t)=\lambda_1e_1(t)+\cdots +\lambda_ne_n(t),
$$
therefore the maps $t\mapsto e_k(t)$, $(1\leq k\leq n)$, are smooth curves 
in the manifolds $\mathfrak{M}(r_k)$ of the projections in $Z$ that 
have fixed finite rank $r_k=\rank (e_k)$, whose tangent spaces at $e_k=e_k(0)$ are  
$Z_{1/2}(e_k)$ (see \cite{CHIS} or \cite{ISMA}). Therefore 
$$u_k\colon = \frac{d}{dt}\vert_{t=0}e_k(t)\in Z_{1/2}(e_k), \qquad 1\leq k\leq n.$$
Since the spectral projections of $a(t)$ corresponding to different 
spectral values $\lambda_k\neq \lambda_j$ are orthogonal, we have $e_j(t)\,e_k(t)=0$ 
for all $t\in I$, and taking the derivative at $t=0$, 
\begin{equation}\label{tv} 
e_j\, u_k=u_k\,ej=0, \qquad j\neq k, \;\; 1\leq j,k\leq n.
\end{equation} 
By \ref{tv}, the tangent vector to $t\mapsto a(t)$ at $t=0$, that is, $u\colon =
\frac{d}{dt}\vert_{t=0}a(t) = \sum_k\lambda_ku_k$ satisfies 
\begin{align}
\{ {\bf a}\, {\bf a}\, u\}=\{ \sum_je_j\, \sum_k e_k\,\sum_l\lambda _lu_l \}&= 
\sum_{j,k,l}\lambda_l\,\{e_je_ku_l\} =\nonumber\\
\sum_{k,l}\lambda_l \{e_ke_ku_l\}= \frac{1}{2}\sum_{k,l}\lambda_l(e_ku_l+u_le_k)&=
\sum_l\lambda_l\{e_le_lu_l\}= \frac{1}{2}\sum_l\lambda_lu_l=\frac{1}{2}u \nonumber
\end{align}
hence $u\in Z_{1/2}({\bf a})$, and $T_aM$ can be identified with a vector 
subspace of $Z_{1/2}({\bf a})$. In fact $T_aM=Z_{1/2}({\bf a})$ as it easily follows 
from the following result that should be 
compared with (\cite{CHIS} th. 3.3)
\begin{Theor}\label{ts}
The sets $M=M(n,\Lambda ,R)$
defined in (\ref{nn}) are holomorphic direct submanifolds of $Z$. 
The tangent space at the point $a\in M$ is the Peirce subspace
$Z_{1/2}({\bf a})$ where ${\bf a}=\supp (a)$, and a local chart at $a$ 
given by 
\begin{equation}\label{lch}
f\colon u \mapsto f(u)\colon = (\exp g({\bf a},u))a
\end{equation}
with $g({\bf a},u)={\bf a}\Box u-u\Box {\bf a}$. 
\end{Theor}
\Dem $M\subset Z$ is invariant under $\Aut^{\circ} (Z)$. 
Fix any $a\in M$ and let 
$X\colon =  \big(\bigoplus_{1\leq k\leq n}i\,A(e_k)\big) \oplus Z_{1/2}({\bf a}).$
Thus $Z=X\oplus Y$ for a certain subspace $Y$. 
The mapping $X\oplus Y\to Z$ defined by $(x,y)\mapsto F(x,y)\colon =
(\exp g({\bf a},x))y\in
Z$ is a real-analytic and its Fr\'echet derivative at $(0,a)$ is 
invertible. In fact this derivative is 
\begin{eqnarray*}
\frac{\partial F}{\partial x}\vert_{ (0,a)}(u,v)&= g(\ba ,u)a=\Phi_a(u), \\
\frac{\partial F}{\partial y}\vert _{(0,a)}(u,v)&= \big(\exp g(\ba ,0)\big) v=v,
\end{eqnarray*}
which is invertible according to (\ref{lh}). 
By the
implicit function theorem there are open sets $U, V$ with $0\in U\subset X$ 
and $a\in V\subset Y$ such that $W\colon = F(U\times V)$ is open in $Z$ and 
$F\colon U\times V\to W$ is bianalytic. 

To simplify notation set $X_1= Z_{1/2}({\bf a})\subset X$. 
Then $f=F\vert X_1$ establishes
a real analytic homeomorphism between the sets $N_1\colon =U\cap X_1$ and 
$M_1\colon = f(N_1)$. Since $X_1$ is a direct summand in $X$ (hence also in
$Z$), the image $M_1=f(N_1)$ is a direct submanifold.

The operator $g({\bf a},x)= {\bf a}\Box x-x\Box
{\bf a}$ is an inner JB$^*$triple derivation of $Z$, hence $h\colon =\exp g({\bf a}, u)$
is a JB$^*$-triple automorphism of $Z$. Actually $h$ lies in $\Aut^{\circ}(Z)$, 
the identity connected component. But it is known (\cite{KAUG}) that $\Aut (Z)$ 
has two connected components and that the elements in the identity component are 
C$^*$-algebra automorphisms of $Z$ since they have the form $z\mapsto UzU^*$ for 
some $U$ in the unitary group of $H$. In particular $h$ preserves normality,   
spectral values and ranks hence it preserves $M$ and so 
$$ M_1=f(N_1)=\{(\exp g({\bf a},u))a \, \colon u\in N_1\}\subset M.$$
To complete the proof, it suffices to show that $f=F\vert X_1$ is a biholomorphic 
mapping. 
The Fr\'echet derivative of $f$ at $a$ is 
$$
f^{\prime}\vert_{a}(u)= g({\bf a} , u)a= \{{\bf a} ,\, u,\,a\}-\{u,\,{\bf a} ,\,a\}, 
\qquad u\in Z_{1/2}(\bf a).
$$
Therefore ${\overline \partial} f^{\prime}u= \{{\bf a} ,\, u,\,a\}$ and 
$\partial f^{\prime}u=-\{u,\,{\bf a} ,\,a\}$ are the (uniquely determined) 
complex-linear and complex-antilinear components of $f^{\prime}u$. The Peirce 
rules give $\{{\bf a} ,\, u\,a\}=0$ for all $u\in Z_{1/2}(\bf a)$, hence 
$f$ is holomorphic and the same argument holds for the inverse $f^{_1}$ map. 
\qed

Remark that if the algebraic element $a$ is a projection then ${\bf a}=a$ and 
$M$ as a differentiable manifold is the one constructed in (\cite{CHIS}
th. 3.3) and \cite{ISMA}.


\section{The Jordan connection on $M(n,\, \Lambda ,\,R)$}

Let $a\in M\colon =M(n,\, \Lambda ,\,R)$ and set ${\bf a}= \supp (a)$. Recall that a
vector field $X$ on $ M$ is a map from $M$ to the tangent bundle $TM$. 
Thus $X _a$, the value of $X$ at $a\in M$, satisfies 
$X _a\in T_a M\approx Z_{1/2}({\bf a})$. We let ${\frak D} (M)$ be the 
Lie algebra of smooth vector fields on $M$. 
Since the tangent space $T_a M$ at $a\in  M$ 
has been identified with $Z_{1/2}({\bf a})$, we shall consider every vector field 
on $M$ as a $Z$-valued function such that the value at $a$ is 
contained in $Z_{1/2}({\bf a})$. 
Let $ Y^{\prime} _a$ 
be the Fr\'echet derivative of $Y \in {\frak D}(M)$ at $a$. 
Thus $ Y^{\prime}_a$ is a bounded linear operator $Z_{1/2}({\bf a}) 
\to Z$, hence 
$ Y^{\prime} _a X_a\in Z$ and it 
makes sense to take the projection $P_{1/2}({\bf a}) Y^{\prime} _a
X_a
\in Z_{1/2}({\bf a})\approx T_a M$. 
\begin{Def}\label{na}
We define a connection $\nabla$ on $ M$ by 
$$(\nabla _ X  Y)_a\colon = P_{1/2}(\supp (a))\,  Y^{\prime}_a
X_a, \qquad  X, \, Y\in {\frak D}(M), \qquad a\in M.$$
\end{Def} 
Note that if $a$ is a projection, then $a= \supp (a)$ and $\nabla$ coincides
with  the affine connection defined in (\cite{CHIS} def 3.6) and \cite{ISMA}.
It is a matter of routine 
to check that $\nabla$ is an affine connection on $M$, that it is $\Aut ^{\circ}(Z)$-
invariant and torsion-free, i. e., 
$$g\,(\nabla _ X Y)= \nabla _{g (X)}\;g\,( Y), \qquad 
g\in \Aut ^{\circ}(Z),$$
where $(g\, X)_a\colon = g^{\prime}_a\,
( X_{g^{-1}_ a})$ for all 
$X \in {\frak D} (M)$, and 
$$ T(X,  Y)\colon = \nabla _ X  Y- \nabla _ Y X- 
[ X Y]=0, \qquad X, Y\in {\frak D}
( M).$$ 
\begin{Theor}\label{geo}
Let the 
manifolds $M$ be defined as in (\ref{nn}). Then the $\nabla$-geodesics of $M$ 
are the curves 
\begin{equation}\label{geo}
\gamma (t)\colon = (\exp \,t\,g({\bf a},u))a, \qquad t\in \mathbb{R},
\end{equation}
where $a\in M$ and $u\in Z_{1/2}(\bf a)$.
\end{Theor}
\Dem 
Recall that the geodesics of $\nabla$ are the curves   
$t\mapsto \gamma (t)\in M$ that satisfy 
the second order ordinary differential equation 
$$\big(\nabla_{\dot \gamma (t)}\,\dot \gamma (t)\big)_{\gamma (t)}=0.$$ 
Let $u\in Z_{1/2}({\bf a})$. Then $g({\bf a}, u)={\bf a}\Box u
-u\Box {\bf a}$ is an inner JB$^*$-triple derivation of $Z$,
and, as established in the proof of (\ref{ts}),  
$h(t)\colon = \exp \,t\,g({\bf a},u)$ 
is an inner C$^*$-automorphism of $Z$. Thus 
$h(t)a\in M$ 
and $t\mapsto \gamma (t)$ is a curve in the manifold $M$. 
Clearly $\gamma (0)=a$ and taking the derivative with respect to $t$ at 
$t=0$ we get by the Peirce
rules
\begin{align}
\dot \gamma (t)&= g({\bf a},u)\gamma (t)= h(t) g({\bf a},u)a, &
\dot \gamma (0)&=g({\bf a},u)a\in Z_{1/2}({\bf a}), \nonumber \\
\ddot \gamma (t)&=g^2 ({\bf a},u)\gamma (t)=h(t)g({\bf a},u)^2a, &
\ddot \gamma (0)&= g({\bf a},u)\dot \gamma (0)\in Z_1({\bf a})\oplus 
Z_0({\bf a}).\nonumber
\end{align}
In particular $P_{1/2}({\bf a})g({\bf a},u)^2a=0$. The definition of $\nabla$ and the
relation (\ref{spi}) give 
\begin{eqnarray*}
\big(\nabla_{\dot \gamma (t)}\,\dot {\gamma} (t)\big)_{\gamma (t)}& = 
P_{1/2}(\supp \gamma (t))\, \Big( \dot {\gamma} (t)^{\prime} _{\gamma (t)}
\,\dot \gamma (t)  \Big)= 
P_{1/2}(\supp \gamma (t))\,\ddot {\gamma} (t)=\\
{}&P_{1/2}\big(\supp h(t)a\big)\, h(t)g({\bf a},u)a= h(t) P_{1/2}\big(\supp
(a)\big)\,g({\bf a},u)^2a=0
\end{eqnarray*}
for all $t\in \mathbb{R}$. Using the representation $u=\sum_ku_k$ given by
(\ref{1/2}) one gets $g({\bf a},u)a=-\frac{1}{2}\sum_k\lambda_ku_k$, and as 
$\lambda\in \sigma (a)\backslash \{0\}$ the mapping $u\mapsto g({\bf a},u)a$ is 
a linear homeomorphism of $Z_{1/2}({\bf a})$. Since
geodesics are uniquely  determined by the initial point
$\gamma (0)=a$ and the initial velocity 
$\dot \gamma (0)=g({\bf a},u)a$, the above shows that family of curves in 
(\ref{geo}) with $a\in M$ and $u\in T_aM\approx Z_{1/2}(\bf a)$ 
are all geodesics of the connection $\nabla$. \qed

Recall that ${\bf a}= \supp (a)$ is a finite rank projection, hence 
by (\cite{ISMA}, th. 1.1) the JB$^*$-subtriple $Z_{1/2}(\bf a)$ has 
finite rank and the tangent space 
$T_aM\approx Z_{1/2}(\bf a)$ is linearly homeomorphic 
to a Hilbert space under an $\Aut ^{\circ} (Z)$-invariant scalar product 
(say $\langle\cdot \, , \cdot \rangle$). Thus we can define a Riemann metric on $M$ by 
\begin{equation}\label{rm}
g_a(X, Y)\colon = \langle X_{\bf a}, Y_{\bf a}\rangle , \qquad X, \, Y\in 
{\mathfrak D}(M), \qquad a\in M.
\end{equation}
Remark that $g$ is {\sl hermitian}, i.e. we have $g_a(iX, \,iY)=g_a(X,Y)$, and that 
it has been defined in algebraic terms, hence it is $\Aut ^{\circ}(Z)$-invariant. 
Moreover, $\nabla$ is compatible with the Riemann structure, i. e. 
$$
X\, g(Y, W)= g(\nabla_X\,Y,\, W)+g(Y,\,\nabla_X\, W), \qquad X,Y,W\in{\mathfrak D}(M).
$$
Therefore, $\nabla$ is the only Levi-Civita connection on $M$. On the other hand, 
let the map $J\colon Z_{1/2}({\bf a})\to Z_{1/2}({\bf a})$ be given by $Jz\colon =i z$. 
Clearly $J^2= -\Id$, hence $J$ defines (the usual) complex structure on the tangent 
space to $M$ and $\nabla$ is $J$-hermitian  
$$
\nabla_X\, (iY)=i\nabla_X\, Y, \qquad X,Y\in {\mathfrak D}(M),
$$
hence $\nabla$ is the only hermitian connection on $M$. Thus the Levi-Civita and 
the hermitian connection are the same in this case, and so $\nabla$ is the K\"ahler 
connection on $M$.

For a tripotent $e\in\Tri (Z)$, the Peirce reflection around $e$ is the  
linear map $S_e\colon = \Id -P_{1/2}(e)$ or in detail $z=z_1+z_{1/2}+z_0\mapsto 
S_e(z) =z_1-z_{1/2}+z_0$ where $z_k$ are the Peirce $e$-projections of $z$,
$(k=1,1/2,0)$. Recall that $S_e$ is an involutory triple automorphism of $Z$ with 
$S_e(e)=e$, and that if $e$ is a projection (taken as a tripotent) then $S_e$ is a  
C$^*$-algebra automorphism of $Z$. This applies to ${\bf a}= \supp (a)$, 
hence to each $a\in M$ we get $S_{\bf a}$, an involutory automorphism of the 
manifold $M$ which in this way becomes a symmetric holomorphic Riemann (K\"ahler) 
manifold. Note that in general ${\bf a}\notin M$ even if $a\in M$, hence $S_{\bf a}$ 
may have no fixed points in $M$. 

It would be interesting to know if any two points $a, b$ in $M$ can be joined 
by a geodesic and whether geodesics are minimizing curves for the Riemann 
distance. The answers to these questions are affirmative when $M$ consists of 
projections of the same finite rank (see \cite{ISMA}).


\section{Algebraic elements in JB$^*$-triples}

The role that projections play in the study of
algebras is taken by tripotents in the study
of triple systems. A spectral calculus and a
notion of algebraic elements is available in
the stetting of JB$^*$-triples. In what
follows we shall consider the manifold of all 
finite rank algebraic elements in a 
JB$^*$-triple $Z$. 
\begin{Def} 
An element $a\in Z$ is called {\sl algebraic}
if there exits a decomposition 
\begin{equation}\label{aet}
a= \lambda_1e_1+\cdots +\lambda_ne_n
\end{equation}
where 
$(e_k)$ is a family of pairwise orthogonal
tripotents in $Z$ and $(\lambda_k)$ are
complex coefficients. 
\end{Def}
For an algebraic element $a\in Z$ the above
decomposition can always be chosen in such a
way that every $e_k$ is non-zero and the
$\lambda_k$ are real numbers with
$0<\lambda_1<\cdots \lambda_n$, and under
these additional conditions the spectral
representation of $a$ is unique. Clearly $a$
has finite rank if and only if every so
does every $e_k$. 

Remark that for $Z= \mathcal{L}(H)$, normal
algebraic elements in the C$^*$-algebra $Z$ are
algebraic elements in $Z$ as a JB$^*$-triple. 
Given a positive integer $n\in \mathbb{N}$, 
an increasing n-uple of non-zero real numbers
$\Lambda =(\lambda_1,\cdots , \lambda_n)$ and 
an n-uple $ R=(r_1,\cdots , r_n)$ where 
$0< r_k\in \mathbb{N}$, we define 
\begin{equation}\label{mm}
N(n,\,\Lambda ,\, R)\colon =
\{\,\sum_k\lambda_ke_k\;\colon \;
 e_j\Box e_k=0 \;\; \hbox{\rm
for}\; j\neq k,\;\; 
\rank (e_k)=r_k, \; 1\leq j,\, k\leq n \,\} 
\end{equation}
to be the set of the elements (\ref{aet}) where
the  coefficients $\lambda_k$ and ranks $r_k$
are given and the $e_k$ range over non-zero, 
pairwise orthogonal tripotents in $Z$ such that 
$\rank (e_k)=r_k$. The set 
$\mathcal{A}$ of finite rank algebraic
elements in $Z$ is the disjoint union 
$\mathcal{A}=\cup_{n,\, \Lambda , \, R}N(n,\,
\Lambda , \,R)$. 
\begin{Lemma}
Let $Z$ be an irreducible JBW$^*$-triple. Then
each of sets $N=N(, n\,\Lambda , \, R)$ is an 
$\Aut^{\circ}(Z)$-invariant connected subset of
$Z$ on which the group
$\Aut^{\circ}(Z)$ acts transitively. 
\end{Lemma}
\Dem Irreducible JBW$^*$-triples are Cartan
factors and we may assume that $Z$ is a not
{\sl special} as otherwise $\dim Z<\infty$ and
the result is known \cite{LOOS}. Thus
$Z$ is  a J$^*$-algebra in the sense of Harris 
\cite{HAR} that is, a weak*-operator closed
complex linear subspace of
$\mathcal{L}(H,K)$ that is closed under the
operation of taking triple products, for
suitable complex Hilbert spaces $H, K$ with
$\dim H\leq \dim K$. Tripotents are the 
partial isometries 
$e\colon H\to K$ that lie in $Z$. 

We make a type by type proof. Let
$Z=\mathcal{L}(H,K)$ be a type I Cartan factor 
and let $a,b\in N$. In particular  
$$
a=\lambda_1e_1+\cdots +\lambda_ne_n, \qquad 
b=\lambda_1e^{\prime}_1+\cdots
+\lambda_ne^{\prime}_n
$$ 
Let $H_k, H_k^{\prime}\subset H$ be the
domains of the partial isometries $e_k$ and
$e^{\prime}_k$, and similarly let
$K_k,K_k^{\prime}\subset K$ denote their
respective ranges. 
Since $e_k$ and $e^{\prime}_k$ have the same
finite rank $r_k$, they are unitarily
equivalent, that is there are unitary 
operators $U_k\colon H_k\to H_k^{\prime}$
and  $V_k\colon K_k\to K_k$ such that 
$e^{\prime}_k = V_ke_kU_k$. Since the $e_k$
are pairwise orthogonal we have $H_k\perp H_j$
and $K_k\perp K_j$ for $k\neq j$ and $\bigoplus
U_k$, $\bigoplus V_k$ are unitary operators on
$\bigoplus H_k$ and $\bigoplus K_k$ that can be
extended to unitary  operators
$U\colon H\to H$ and $V\colon K\to K$ if
needed. The mapping $Z\to Z$ given by
$z\mapsto VzU$ is a JB$^*$-triple automorphism
that lies in $\Aut^{\circ}(Z)$ \cite{KAUG}
and clearly satisfies $b=VaU$. Hence
$\Aut^{\circ}(Z)$ acts transitively on $N$, 
$N$ is connected and invariant under that
group. 

Cartan factors of types II and III can
treated in the same way. The case of spin
factors may be discussed with a different
approach, but we shall not go into details. 
\qed

Now consider the joint Peirce decomposition of $Z$ relative to the 
family $(e_1, \cdots , e_n)$ where $a= \lambda_1e_1+ \cdots +\lambda_ne_n$ 
is the spectral resolution of $a$. Let the {\sl support} of $a$ 
be tripotent ${\bf a}= \supp a\colon = e_1+\cdots
+e_n$, and  note that 
$$X\colon =  \big(\bigoplus_{1\leq k\leq n}i\,A(e_k)\big) \oplus 
Z_{1/2}({\bf a}).$$
is a topologically complemented subspace in $Z$.

Fix one of the sets $N=N(n,\Lambda ,R)$ and a point $a\in N$ with spectral 
resolution $a= \sum_k\lambda_ke_k$. From the properties $e_k\Box e_j=0$ 
for $j\neq k$, the successive odd powers of $a$ have the
expression   
$$ a^l= \lambda_1^{2l+1}e_1+\cdots + \lambda_n^{2l+1}e_n, \qquad 0\leq
l\leq n-1,$$ where the determinant $\det (\lambda_k^{2l+1})\neq 0$ does
not vanish since it is a Vandermonde determinant and the $\lambda_k$ are
pairwise distinct. Thus the $e_k$ are polynomials in $a$ whose   
coefficients are rational functions of the $\lambda_k$.  
Suppose $N$ is a differentiable manifold, and let us obtain its tangent
space
$T_aN$. Consider a smooth curve $t\mapsto a(t)$ in $N$ through $a$, $t\in
I$, for a neighbourhood $I$ of $0\in \mathbb{R}$ and $a(0)=a$. Each
$a(t)$ has a spectral resolution 
$$
a(t)=\lambda_1e_1(t)+\cdots +\lambda_ne_n(t),
$$
therefore the maps $t\mapsto e_k(t)$, $(1\leq k\leq n)$, are smooth 
curves in the manifolds $\mathfrak{N}(r_k)$ of the tripotents
in $Z$ that have fixed finite rank $r_k=\rank (e_k)$, whose tangent
spaces at
$e_k=e_k(0)$ are respectively 
$i\,A(e_k)\oplus Z_{1/2}(e_k)$ (see \cite{CHIS} or \cite{ISMA}). Therefore 
$$ z_k\colon =\frac{d}{dt}\vert_{t=0}\,e_k(t)= iv_k+u_k\colon \in
i\,A(e_k)\oplus Z_{1/2}(e_k), \qquad 1\leq k\leq n.$$ 
Set $v\colon = \sum_k \lambda_kv_k$ and $u\colon = \sum_k\lambda_ku_k$. From 
$Z_1(e_k)\Box Z_0(e_j)=\{0\}$, we get 
$$ 
\{ {\bf a}\, {\bf a}\, iv\}=i \sum_{j,k,l}\lambda_l\{e_je_kv_l\}=
i\sum_k\lambda_kv_k=iv\in i\, \bigoplus_k A(e_k)
$$ 
The spectral
tripotents of $a(t)$ corresponding to different spectral values 
$\lambda_k\neq \lambda_j$ are orthogonal, hence $e_j(t)\Box e_k(t)=0$ 
for all $t\in I$, and taking the derivative at $t=0$ we get 
\begin{equation}\label{tv} 
e_j\Box u_k=u_k\Box ej=0, \qquad j\neq k, \;\; 1\leq j,k\leq n.
\end{equation} 
Hence 
$$
\{ {\bf a}\, {\bf a}\, u\}=\{ \sum_je_j\, \sum_k e_k\,\sum_l\lambda_lu_l\}= 
\sum_{j,k,l}\lambda_l\{e_je_ku_l\}=\frac{1}{2} \sum_k\lambda_ku_k=\frac{1}{2}u
$$
which shows that $u\in Z_{1/2}({\bf a})$. 
By \ref{tv}, the tangent vector to $t\mapsto a(t)$ at $t=0$ is  
$z\colon =
\frac{d}{dt}\vert_{t=0}a(t) = \sum_k\lambda_k(iv_k+u_k)=iv+u$ hence it satisfies 
$$
\{ {\bf a}\, {\bf a}\, z\}=iv+\frac{1}{2}u\in 
i\bigoplus A(e_k)\oplus  Z_{1/2}({\bf a}),$$ 
hence $T_aN$ can be identified with a vector 
subspace of $i\bigoplus A(e_k)\oplus  Z_{1/2}({\bf a})$. In fact $T_aN$ 
coincides with that space as it easily follows from the following result that
should be  compared with (\cite{CHIS} th. 3.3)

\begin{Theor}\label{tst}
The sets $N=N(n,\Lambda ,R)$
defined in (\ref{mm}) are real analytic direct submanifolds of $Z$. 
The tangent space at the point $a\in N$ is the Peirce subspace
$X$, where
${\bf a}=\supp (a)$, and a local chart at $a$  given by 
\begin{equation}\label{lcht}
f\colon z \mapsto f(z)\colon = (\exp g({\bf a},z))a
\end{equation}
with $g({\bf a},z)={\bf a}\Box z-z\Box {\bf a}$. 
\end{Theor}
\Dem $N\subset Z$ is invariant under $\Aut^{\circ} (Z)$. 
Fix any $a\in N$ and let 
$X\colon =  \big(\bigoplus_{1\leq k\leq n}i\,A(e_k)\big) \oplus Z_{1/2}({\bf a}).$
Thus $Z=X\oplus Y$ for a certain subspace $Y$. 
The mapping $X\oplus Y\to Z$ defined by $(x,y)\mapsto F(x,y)\colon =
(\exp g({\bf a},x))y\in
Z$ is a real-analytic and its Fr\'echet derivative at $(0,a)$ is 
invertible as proved in (\ref{lh}). 
By the
implicit function theorem there are open sets $U, V$ with $0\in U\subset X$ 
and $a\in V\subset Y$ such that $W\colon = F(U\times V)$ is open in $Z$ and 
$F\colon U\times V\to W$ is bianalytic and the image $F(U)$ is a direct 
real analytic submanifold of $Z$.

The operator $g({\bf a},z)= {\bf a}\Box z-z\Box
{\bf a}$ is an inner JB$^*$triple derivation of $Z$, hence 
$h\colon =\exp g({\bf a}, z)$
is a JB$^*$-triple automorphism of $Z$. Actually $h$ lies in $\Aut^{\circ}(Z)$, 
the identity connected component. In particular $h$ preserves the algebraic   
character and the spectral decomposition, hence it preserves $N$ and so 
$$ F(N)=\{(\exp g({\bf a},z))a \, \colon z\in U\}\subset N.$$

This completes the proof. \qed

\begin{Def} For the tripotents $e, \, e^{\prime}$ we set 
$e\sim e^{\prime}$ if and only if $e$ and $e^{\prime}$ have the same 
$k$-Peirce projectors for $k=0,1/2,1$. 
\end{Def}

This notion was introduced by Neher who proved 
(\cite{NEH}, th.2.3) that
\begin{equation}\label{eq} 
e\sim e^{\prime}\;\; \Longleftrightarrow \;\; 
e\in Z_1(e^{\prime}) \;\;\hbox{\rm and}\;\; 
e^{\prime}\in Z_1(e),
\end{equation} 
or equivalently if and only if 
$e\square e=
e^{\prime}\square e^{\prime}$. Next we extend this relation to an 
equivalence in the manifold $N$. 

\begin{Def} Let $a, b$ be elements in $N$ with spectral resolutions 
$a=\sum_k\lambda_ke_k$ and $b=\sum_k\lambda_kf_k$ respectively. We say 
that $a$ and $b$ are equivalent (and write $a\sim b$) if the joint Peirce 
decompositions of $Z$ relative to the orthogonal families $\mathcal{E} = 
(e_k)$ and $\mathcal{F}=(f_k)$ are the same. 
\end{Def}
Note that $\sim$ coincides with the equivalence of Neher when the algebraic
elements $a$ and $b$ are tripotents. By (\cite{LOOS}, th. 3.14), the 
Peirce spaces of the tripotent $e_k$ can be expressed in terms of the 
joint Peirce decomposition of
$Z$  relative to $\mathcal{E}$, hence $a\sim b$ if and only if $e_k\sim f_k$ 
for $1\leq k\leq n$. 

\begin{Prop}\label{fb}
Let $a,b$ be points in $N$ such that $a=\sum\lambda_ke_k$ and 
$b=(\exp g({\bf a}, z)a$ for some tangent vector $z= iv+u\in 
\big(\bigoplus_{1\leq k\leq n}i\,A(e_k)\big) \oplus Z_{1/2}({\bf a})$. 
Then $a\sim b$ if and only if $u=0$.  
\end{Prop}
\Dem Let $b=(\exp g({\bf a}, z)a=\sum_k\lambda_kf_k$ be the spectral 
resolution of $b$. Then each $f_k$ is an odd polynomial in $b$, say 
$f_k=p_k(b)$, $1\leq k\leq n$. To simplify the notation, consider the index $k=1$ and omit 
the reference to it in the rest of the proof. If $a\sim b$ then $e\sim f$ 
hence by (\ref{eq}) we must have $f=\{eef\}$ that is 
\begin{equation}\label{ecc}
p(b)=\{eep(b)\}=p(\{eeb\})
\end{equation} 
Clearly we have $\rho b\sim a$ for all $\rho\in \mathbb{T}$, which replaced
above yields an identity between two polynomials in $\rho$. 
Let $X^m$, for some positive odd integer $m$, be the term of $p$ of lowest 
degree whose coefficient is not zero. Then (\ref{ecc}) entails 
$b^m=\{eeb^m\}$, that is 
$(\exp g({\bf a}, z))^m a= \{e\,e\,(\exp g({\bf a}, z))^m a \}.$ 
Taking the Fr\'echet derivative at the origin 
$g({\bf a}, \cdot )\,a =\{e\,e\,g({\bf a}, \cdot )\,a\},$ 
which evaluated at the tangent vector $z=iv+u=i\sum_kv_k+\sum_ku_k$ and 
using the Peirce rules as in the proof of (\ref{lh}) yields $u=0$. The 
converse is easy. \qed 

In particular, there is a neighbourhood of $a$ in $N$ in which the 
algebraic elements $b$ equivalent to $a$ are those of the form $b=( \exp 
g({\bf a},\, iv))\,a$ with $v=\sum_kv_k\in \bigoplus_kA(e_k)$, 
which gives the expression of the fibre of $N$ through $a$. 
\begin{Prop} 
Let $a\in N$ be an algebraic element in $Z$ with spectral resolution 
$a=\sum_k\lambda_ke_k$. Then the fibre of $N$ through $a$ is the set 
of the elements $\sum_k\lambda_k z_k$ where $z_k$ lies in the unit circle 
of the JB$^*$-algebra $Z_1(e_k)$ for $1\leq k\leq n$.
\end{Prop}
\Dem Let $v=\sum_k\lambda_kv_k\in \bigoplus_kA(e_k)$, and consider the curves 
in $Z$ 
$$
\phi (t) \colon = (\exp \,t g({\bf a}, iv))a, \qquad 
\psi(t) \colon = \sum_k\,\lambda_k(\exp \, t g(e_k, iv_k))e_k\colon = \sum_k\lambda_k\psi_k(t), 
\qquad t\in \mathbb{R}.
$$
They are the solutions of the differential equations 
$$
\frac{d\phi (t)}{dt}= g({\bf a}, \phi (t)), \qquad 
\frac{d\psi (t)}{dt}= \sum_k\,\lambda_kg(e_k, \psi_k (t))
$$ 
with the initial conditions $\phi (0)=a$ and
$\psi(0)=\sum_k\lambda_ke_k=a$ respectively.  From $Z_1(e_k)\Box Z_1(e_j)=\{0\}$ for $k\neq j$ we 
get 
$$
g({\bf a}, iv)= g(\sum_ke_k, \, i\sum_j \lambda_jv_j)= \sum_k\lambda_kg(e_k,\,iv_k) 
$$
and the uniqueness of solutions of differential equations gives 
$\phi (t) = \sum_k\lambda_k \psi_k (t)$ for all $t\in \mathbb{R}$.  But it is 
known (\cite{LOOS} th. 5.6) that for fixed $k$, $1\leq k\leq n$, the set 
$z_k=(\exp t g(e_k, iv_k))e_k$, $t\in \mathbb{R}$, $v_k\in A(e_k)$, is the unit circle of the
JB$^*$-algebra 
$Z_1(e_k)$, that is the set of those $w\in Z_1(e_k)$ that satisfy 
$w^* =w^{-1}$. This completes the proof.\qed   

By restricting the local charts in (\ref{lcht}) to the direct summand 
$Z_{1/2}({\bf a})\subset T_aN$ we get a direct submanifold $B=B(n, \Lambda , R)$ of 
$Z$, and we refer to $B$ as the {\sl base} manifold of $N$. Clearly $B$ is 
a holomorphic submanifold of the real analytic manifold $N$, and as in 
section 3 
$$
\big(\nabla_X\, Y \big)_a\colon = P_{1/2}({\bf a}) Y^{\prime}_aX_a, \qquad 
X, Y\in \mathfrak{D}(B), \quad a\in B,
$$
is an $\Aut^{\circ}(Z)$-invariant torsionfree affine connection on $B$ 
whose geodesics are the curves $\gamma(t)\colon = (\exp t\,g({\bf a}, u))a$, 
$t\in \mathbb{R}$, for $a\in B$ and $u\in Z_{1/2}({\bf a})$. Moreover, for 
$a\in B$ the Peirce reflection with respect to ${\bf a}$ is an involutory triple  
automorphisms of $Z$ that fixes ${\bf a}$, hence it fixes $i\, \bigoplus_k 
A(e_k)$ and $Z_{1/2}({\bf a})$. It is easy to see that this reflection 
commutes with the exponential mapping, hence it fixes $B(n, \Lambda , R)$ 
and os it defines a holomorphic symmetry of $B$. In general $({\bf a})$ 
does not belong to $B$ hence this symmetry in general has no fixed points in $B$. 
When the algebraic element $a\in Z$ has finite rank, that is when $\rank (a) = 
\sum_k\rank (e_k)<\infty$, the subtriple $Z_{1/2}({\bf a})$ is linearly equivalent 
to a complex Hilbert space by \cite{KAUS} and by using the algebraic metric of 
Harris one can introduce an $\Aut^{\circ}(Z)$-invariant Riemann structure and a 
K\"ahler structure on the base manifold in exactly the same way we did in section 
3, and the connection $\nabla$ turns out to be the Levi-Civita and the K\"ahler 
connection on $B$.


\end{document}